\documentclass[a4paper,12 pt]{amsart}

\usepackage{amsfonts}
\usepackage{amsmath}
\usepackage{amsmath,amsfonts,amssymb,amsthm,latexsym}
\usepackage{graphicx}
\usepackage{hyperref}
\usepackage{systeme}

\newcommand{\RR}{\mathbb{R}}
\newcommand{\CC}{\mathbb{C}}
\newcommand{\NN}{\mathbb{N}}

\newcommand{\RE} { {\rm Re \,} }
\newcommand{\IM} { {\rm Im \,} }

\newtheorem{Tw}{Theorem}
\newtheorem{Le}{Lemma}
\newtheorem{Stw}{Proposition}

\newtheorem{Wn}{Corollary}
\DeclareMathOperator{\lin}{lin}

\theoremstyle{remark}
\newtheorem{Uw}{Remark}

\theoremstyle{definition}
\newtheorem{Df}{Definition}

\begin{document} 
\title[On polyharmonic polynomials]{On  polyharmonic polynomials}
\author{Hubert Grzebu{\l}a}
\address{Faculty of Mathematics and Natural Sciences,
College of Science\\
Cardinal Stefan Wyszy\'nski University\\
W\'oycickiego 1/3,
01-938 Warszawa, Poland}
\email{h.grzebula@uksw.edu.pl}
\author{S{\l}awomir Michalik}
\address{Faculty of Mathematics and Natural Sciences,
College of Science\\
Cardinal Stefan Wyszy\'nski University\\
W\'oycickiego 1/3,
01-938 Warszawa, Poland}
\email{s.michalik@uksw.edu.pl}
\urladdr{\url{http://www.impan.pl/~slawek}}

\begin{abstract}
We study the orthogonal projection of homogeneous polynomials onto the space of homogeneous polyharmonic polynomials. To do this we derive the decomposition of homogeneous polynomials in terms of the Kelvin transform of derivatives of the fundamental solution  $|x|^{2-n}$ or $\log |x|$. We consider also the vector bases of the space of homogeneous polyharmonic polynomials and study the problem of convergence  of orthogonal series.
\end{abstract}
\maketitle   
\section{Introduction}
 The theory of polyharmonic functions has been investigated extensively. The basic work about this class of functions is due to Aronszajn, Creese and Lipkin \cite{A-C-L}, and it is still an active field of interest in mathematics, see for example \cite{G-G-S} and references therein. These functions  are also very useful in many branches in applied mathematics, such as approximation theory, polysplines, radial basis functions and wavelet analysis 
(see for example \cite{Bu,K,M-N}). 
 
In this paper we study the polyharmonic polynomials. A lot of information about them can be found in a classical paper \cite{B-C}. We start our considerations with the decomposition of homogeneous polynomials. There are many works about the decomposition of the homogeneous polynomials. In the papers \cite{L-R,R3} the Fischer decomposition is considered, and it is relevant to the well known Khavinson-Shapiro conjecture (see also \cite{R1,R2}). The  decomposition given in this paper is a natural extension of the result from \cite{A-B-R} (see also \cite{A-R}). It is worth to mention that the decomposition given there  allows to obtain the solution of many Dirichlet-type problems with polynomial data  without using integration, namely mentioned results lead to the algorithm which involves only differentiation (see \cite{A-R}).

Let $\mathcal{P}_m(\RR^n)$ denote the space of the homogeneous polynomials of degree $m$ on $\RR^n$ with $n>2$. It is known that for any $u\in\mathcal{P}_m(\RR^n)$ the following formula holds (see Lemma 5.17 in \cite{A-B-R})
$$K[u(D)|x|^{2-n}]=c_mu+|x|^2v,$$
where 
$$c_m=\prod_{k=0}^{m-1}(2-n-2k),$$
 $u(D)=\sum_\alpha a_\alpha D^\alpha$ for $u=\sum_\alpha a_\alpha x^\alpha$, and $K$ is the Kelvin transform. In this paper we derive an explicit formula for the polynomial $v\in\mathcal{P}_{m-2}(\RR^n)$:
$$v=\sum_{k=1}^{\left[\frac{m}{2}\right]}\sum_{\substack{(k_1,\cdots,k_n)\in\NN^n_0\\k_1+\dots+k_n=k}}\frac{c_{m-k}}{2^kk_1!\dots k_n!}|x|^{2(k-1)} D^{2k}_{x_1^{2k_1},\dots,x_n^{2k_n}}u.$$
As a corollary we find the mentioned decomposition of the homogeneous polynomials in terms of the Kelvin transform, namely
$$u(x)=\sum_{k=0}^{\left[\frac{m}{2}\right]}\frac{|x|^{2k}}{c_{m-2k}}K[u_k(D)|x|^{2-n}],$$ 
where $u_0=u$ and 
$$
u_k=\sum_{\substack{ (k_1,\dots,k_n)\in\NN_0^n\\k_1+\dots +k_n=k}}d_{m,k,k_1,\dots,k_n}D^{2k}_{x_{1}^{2k_1}\dots x_{n}^{2k_n}}u\quad\textrm{for}\quad k=1,\dots,[m/2]
$$
with $d_{m,k,k_1,\dots,k_n}$ being some constants independent of $u$. Analogous results are obtained for $n=2$ by changing the fundamental solutions, here we replace  $|x|^{2-n}$ by $\log|x|$. From the above we easily find the orthogonal projection of homogeneous polynomials onto the space of polyharmonic polynomials and the orthogonal projection of homogeneous polynomials on $\CC^n$  onto the space of so called spherical polyharmonics $\mathcal{H}^p_m(\widehat{S}_p)$. The spherical polyharmonics are a natural generalisation of the well known spherical harmonics. Here we restrict the polyharmonic polynomials to the union of rotated spheres (see \cite{G-M-2}):
$$\widehat{S}_p:=\bigcup_{k=0}^{p-1}e^{\frac{k\pi i}{p}}S.$$ 
The motivation for the study of spherical polyharmonics on the union of rotated balls comes from the Pizzetti formulas for  the polyharmonic operator $\Delta^p$ given in \cite{L0,S-M} (see also \cite{G-M}).

In the next part of this paper we deal with the vector bases of the space  of the polyharmonic homogeneous polynomials of order $p$ and of degree $m\leq p$ on $\CC^n$.  

A criterion for the convergence of some  orthogonal series is considered  in \cite{F-V}. An explicit formula for the reproducing kernel of the Hilbert space $\mathcal{H}_m^p(\RR^n)$ of homogeneous polyharmonic polynomials endowed  with the Fischer inner product is derived in \cite{R}. By this formula, the condition given in \cite{F-V} is improved in \cite{R}. In our paper we study a similar problem. Using \cite{G-M} we get the analogous result for the space of polyharmonic polynomials, but endowed with the integral inner product induced from $L^2(\widehat{S}_p)$. In this paper we also discuss the connection between 
these two kinds of inner products.

\section{Preliminaries}
In this part we recall some notations, definitions and facts that we will use in this paper. By
\begin{equation*}
|x|= \left( \sum_{j=1}^nx_j^2 \right)^{1/2}
\end{equation*}
we denote the real norm of $x=(x_1,\dots,x_n)\in \RR^n$. We also will use the real norm for complex vectors $z=(z_1,\dots,z_n)\in\CC^n$:
\begin{gather}
\label{eq:z}
|z|=\left( \sum_{j=1}^nz_j^2 \right)^{1/2},
\end{gather}
here by a square root we mean the principal square root, where a branch cut is taken along the non-positive real axis.

We will also use the complex norm of $z=(z_1,\dots,z_n)\in\CC^n$ defined by
$$|z|_{\CC^n}=\left(\sum_{j=1}^nz_j\overline{z_j} \right)^{1/2}.
$$

Let $B$ and $S$ be respectively the unit ball and sphere in $\RR^n$ with a centre at the origin. For  the angle $\varphi\in\RR$ we will consider a rotated unit  ball  in $\RR^n$ defined by
$$e^{i\varphi}B:=\{e^{i\varphi}x:x\in B\}.$$
We define similarly a rotated sphere $e^{i\varphi}S$. Let $p\in\NN$. We will denote unions of rotated balls and spheres  as follows
$$
\widehat{B}_p:=\bigcup_{l=0}^{p-1}e^{\frac{l\pi i}{p}}B\quad\textrm{and}\quad\widehat{S}_p:=\bigcup_{l=0}^{p-1}e^{\frac{l\pi i}{p}}S.
$$
More generally, for any $r>0$ we use the notation $\widehat{B}_p(r):=\bigcup_{l=0}^{p-1}e^{\frac{l\pi i}{p}}B(r)$, where $B(r)$ denotes the ball in $\RR^n$ with a radius $r$ and a centre at the origin.

Let $G\subset \RR^n$ be open. A function $u:G\rightarrow\CC$ is called \emph{polyharmonic} of order $p\in\NN$ if $u\in\mathcal{C}^{2p}(G)$ and $\Delta^p u(x)=0$ for $x\in G$, where
$$\Delta=\frac{\partial^2}{\partial x_1^2}+\dots+\frac{\partial^2}{\partial x_n^2}$$
is the \textit{Laplace operator} and $\Delta^p$ is its $p$-th iterate. A polyharmonic function of order $1$ is called \textit{harmonic}. It is well known that
if $u$ is polyharmonic on $B$, then $u$ can be holomorphically extended to the Lie ball (\cite[Theorem D]{S}), in particular $u$ can be extended to  $e^{i\varphi}B$ for every angle $\varphi\in\RR$ (see also \cite[Lemma 1]{G-M}). 

\begin{Le}[Almansi expansion, {\cite[Proposition 1.2, Proposition 1.3]{A-C-L}}]
\label{Almansi}
A function $u$ is polyharmonic on $B$ if and only if there exist unique harmonic functions  $h_0,h_1,\dots, h_{p-1}$ on $B$ such that
\begin{gather*}
u(x)=h_0(x)+|x|^2h_1(x)+\dots+|x|^{2(p-1)}h_{p-1}(x) \quad \textrm{for} \quad x\in B.
\end{gather*}
\end{Le}

\begin{Df}
 A map 
 $$x\mapsto x^*:=\left\{\begin{array}{lll}
  x/|x|^2 &\textrm{if}& x\neq 0,\infty\\
  0 &\textrm{if}& x=\infty\\
  \infty &\textrm{if}& x=0
  \end{array}
  \right.
$$
is called the \emph{inversion} of $\RR^n\cup\{\infty\}$ relative to the unit sphere.
 
 Let $u$ be a function defined on a set $E\subseteq(\RR^n\cup\{\infty\})\setminus\{0\}$. A function $K[u]$ defined on $E^*:=\{x^*\colon x\in E\}$ by
 \begin{equation*}
K[u](x):=|x|^{2-n}u\left(x^*\right)
\end{equation*}
is called the \emph{Kelvin transform} of $u$.
\end{Df}
\begin{Uw}
 We can extend the definition of the Kelvin transform to the complex case when $E\subseteq(\CC^n\cup\{\infty\})\setminus\{z\in\CC^n\colon |z|=0\}$.
 
 Observe that the set $\{z\in\CC:|z|=0\}$  contains not only the origin since $|z|$ given by (\ref{eq:z}) is a complex valued entity. 
\end{Uw}

Let us recall a definition of harmonicity at $\infty$.

\begin{Df}
Let $E\subset\RR^n$ be compact. If $u$ is harmonic on $\RR^n\setminus E$, then $u$ is \textit{harmonic at $\infty$} provided $K[u]$ has a removable singularity at the origin. 
\end{Df}

\begin{Uw}[{\cite[Theorem 4.7, Theorem 4.9]{A-B-R}}]
Let $E\subset\RR^n$ be compact and $u$ be harmonic on $\RR^n\setminus E$. Then $u$ is harmonic at $\infty$ if and only if 
\begin{enumerate}
\item[(i)] $\lim_{x\rightarrow\infty}u(x)$ is finite, when $n=2$,
\item[(ii)] $\lim_{x\rightarrow\infty}u(x)=0$, when $n>3$.
\end{enumerate}
\end{Uw}

\begin{Le}[{\cite[Proposition 1.4]{A-C-L}}]
If $u$ is a harmonic function on $B$, then its Kelvin transform $K[u]$ is harmonic on $(\RR^n\setminus B)\cup\{\infty\}$.
\end{Le}

Let $m,p\in \NN$. We denote by  $\mathcal{P}_m(\RR^n)$ (resp. $\mathcal{P}_m(\CC^n)$) the space of homogeneous polynomials on $\RR^n$ (resp. $\CC^n$) of degree $m$. By $\mathcal{H}_m^p(\CC^n)$ we mean the space of polynomials on $\CC^n$, which are homogeneous of degree $m$ and are polyharmonic of order~$p$. It is easy to see that if $m<2p,$ then $ \mathcal{H}_m^p(\RR^n)=\mathcal{P}_m(\RR^n)$ (resp. $ \mathcal{H}_m^p(\CC^n)=\mathcal{P}_m(\CC^n)$).

\section{Projection of the homogeneous polynomials onto the polyharmonic polynomials}
In this section we prove the main result on the orthogonal projection of the homogeneous polynomials onto the space of polyharmonic polynomials.
The result is based on a more precise version of \cite[Lemma 5.17]{A-B-R}. 
\begin{Stw}
\label{prop:1}
Let $u\in\mathcal{P}_m(\RR^n)$ with $m\geq 1.$

\begin{enumerate}
\item[(i)] If $n\ne 2$, then
$$K[u(D)|x|^{2-n}]=\sum_{k=0}^{\left[\frac{m}{2}\right]}\sum_{\substack{(k_1,\cdots,k_n)\in\NN^n_0\\k_1+\dots+k_n=k}}\frac{c_{m-k}}{2^kk_1!\dots k_n!}|x|^{2k} D^{2k}_{x_1^{2k_1},\dots,x_n^{2k_n}}u,$$
where
$$c_{m}=\prod_{k=0}^{m-1}(2-n-2k).$$
\item[(ii)] If $n=2$, then
$$K[u(D)\log|x|]=\sum_{k=0}^{\left[\frac{m}{2}\right]}\sum_{\substack{(k_1,k_2)\in\NN^2_0\\k_1+k_2=k}}\frac{\tilde{c}_{m-k}}{2^kk_1! k_2!}|x|^{2k} D^{2k}_{x_1^{2k_1},x_2^{2k_2}}u,$$
where
$$\tilde{c}_m=(-2)^{m-1}(m-1)!.$$
\end{enumerate}
\end{Stw}
\begin{proof}
We shall prove (i). Since the Kelvin transform is linear, it is enough to prove it for a monomial $u\in\mathcal{P}_m(\RR^n)$. Let $u(x)=x^\alpha$, $|\alpha|=m$. We proceed by induction on $m$. If $m=1$, then $u(x)=x_i$ for some $i=1,\dots,n$. So $u[D]=D_{x_i}$,
$$K[u(D)|x|^{2-n}]=K[(2-n)|x|^{-n}x_i]=(2-n)|x|^{2-n}\left|\frac{x}{|x|^2}\right|^{-n}\cdot\frac{x_i}{|x|^2}=(2-n)x_i$$
and $c_{1}=2-n$. Hence the proposition is valid.
Let us assume that the proposition holds for $m\geq 1$. Let $u(x)=x^\alpha$, $|\alpha|=m$. By the induction hypothesis we have
$$K[D^\alpha|x|^{2-n}]=\sum_{k=0}^{\left[\frac{m}{2}\right]}\sum_{\substack{(k_1,\cdots,k_n)\in\NN^n_0\\k_1+\dots+k_n=k}}\frac{c_{m-k}}{2^kk_1!\dots k_n!}|x|^{2k} D^{2k}_{x_1^{2k_1},\dots,x_n^{2k_n}}x^\alpha.$$
Since $u\in\mathcal{P}_m(\RR^n),$ so $D^{2k}_{x_{i_1}^2,\dots,x_{i_k}^2}x^\alpha\in \mathcal{P}_{m-2k}(\RR^n)$. Moreover the Kelvin transform is an involution and is linear, therefore using this transform to the last equation we get
\begin{eqnarray*}
D^\alpha|x|^{2-n}&=& \sum_{k=0}^{\left[\frac{m}{2}\right]}\sum_{\substack{(k_1,\cdots,k_n)\in\NN^n_0\\k_1+\dots+k_n=k}}\frac{c_{m-k}}{2^kk_1!\dots k_n!}K[|x|^{2k} D^{2k}_{x_1^{2k_1},\dots,x_n^{2k_n}}x^\alpha]\\
&=& \sum_{k=0}^{\left[\frac{m}{2}\right]}\sum_{\substack{(k_1,\cdots,k_n)\in\NN^n_0\\k_1+\dots+k_n=k}}\frac{c_{m-k}}{2^kk_1!\dots k_n!}|x|^{2-n-2m+2k} D^{2k}_{x_1^{2k_1},\dots,x_n^{2k_n}}x^\alpha
\end{eqnarray*}
Set $\tilde{x}:=(x_2,\dots,x_n)$, $\tilde{\alpha}:=(\alpha_2,\dots,\alpha_n)$. Without loss of generality we may differentiate the last equality with respect to $x_1$:
\begin{multline*}
D_{x_1}D^\alpha|x|^{2-n}=\sum_{k=0}^{\left[\frac{m}{2}\right]}\sum_{\substack{(k_1,\cdots,k_n)\in\NN^n_0\\k_1+\dots+k_n=k}}\frac{c_{m-k}}{2^kk_1!\dots k_n!}\\
\times\left((2-n-2m+2k)|x|^{2-n-2(m+1)+2k} x_1D^{2k}_{x_1^{2k_1},\dots,x_n^{2k_n}}x^\alpha+|x|^{2-n-2m+2k}D^{2k+1}_{x_1^{2k_1+1},\dots,x_n^{2k_n}}x^\alpha\right)\\
=c_{m}(2-n-2m)|x|^{2-n-2m-2}x_1 x^\alpha+I_1+I_2,
\end{multline*}
where
\begin{multline*}
I_1=\sum_{k=1}^{\left[\frac{m}{2}\right]}\sum_{\substack{(k_1,\cdots,k_n)\in\NN^n_0\\k_1+\dots+k_n=k}}\frac{c_{m-k}}{2^kk_1!\dots k_n!}\\
\times(2-n-2m+2k)|x|^{2-n-2(m+1)+2k} x_1D^{2k}_{x_1^{2k_1},\dots,x_n^{2k_n}}x^\alpha
\end{multline*}
and
\begin{equation*}
I_2=\sum_{k=0}^{\left[\frac{m}{2}\right]}\sum_{\substack{(k_1,\cdots,k_n)\in\NN^n_0\\k_1+\dots+k_n=k}}\frac{c_{m-k}}{2^kk_1!\dots k_n!}|x|^{2-n-2m+2k}D^{2k+1}_{x_1^{2k_1+1},\dots,x_n^{2k_n}}x^\alpha.
\end{equation*}
Let us note that
\begin{equation}
\label{f1}
c_{m-k}(2-n-2m+2k)=c_{m+1-k}
\end{equation}
for every $k=0,1,\dots, \left[\frac{m}{2}\right]$. Since $x^\alpha=x_1^{\alpha_1}\tilde{x}^{\tilde{\alpha}}$, we have
\begin{eqnarray}
\label{e4}
\nonumber x_1 D^{2k}_{x_1^{2k_1},\dots,x_n^{2k_n}}x^\alpha &=& x_1 D^{2k_1}_{x_{1}^{2k_1}}x_1^{\alpha_1} D^{2(k-k_1)}_{x_2^{2k_2},\dots,x_n^{2k_n}}\tilde{x}^{\tilde{\alpha}}\\
\nonumber &=& \frac{\alpha_1-2k_1+1}{\alpha_1+1}D^{2k_1}_{x_{1}^{2k_1}} D^{2(k-k_1)}_{x_2^{2k_2},\dots,x_n^{2k_n}}x_1^{\alpha_1+1}\tilde{x}^{\tilde{\alpha}}\\
&=& \frac{\alpha_1-2k_1+1}{\alpha_1+1} D^{2k}_{x_1^{2k_1},\dots,x_n^{2k_n}}x_1 x^{\alpha}.
\end{eqnarray}
By (\ref{f1}) and (\ref{e4}) we get 
\begin{multline}
\label{eI1}
I_1=\sum_{k=1}^{\left[\frac{m}{2}\right]}\sum_{\substack{(k_1,\cdots,k_n)\in\NN^n_0\\k_1+\dots+k_n=k}}\frac{c_{m+1-k}}{2^kk_1!\dots k_n!}\\
\times|x|^{2-n-2(m+1)+2k} \frac{\alpha_1-2k_1+1}{\alpha_1+1} D^{2k}_{x_1^{2k_1},\dots,x_n^{2k_n}}x_1 x^{\alpha}.
\end{multline}
Now we can consider $I_2$. Let us note that
$$D^{2k+1}_{x_1^{2(k_1+1)},\dots,x_n^{2k_n}}x^\alpha=\frac{1}{\alpha_1+1}D^{2(k+1)}_{x_1^{2k_1},\dots,x_n^{2k_n}}x_1x^\alpha$$
so
$$I_2=\sum_{k=0}^{\left[\frac{m}{2}\right]}\sum_{\substack{(k_1,\cdots,k_n)\in\NN^n_0\\k_1+\dots+k_n=k}}\frac{c_{m-k}}{2^kk_1!\dots k_n!}|x|^{2-n-2m+2k}\frac{1}{\alpha_1+1}D^{2(k+1)}_{x_1^{2(k_1+1)},\dots,x_n^{2k_n}}x_1x^\alpha .$$
Replacing $k_1$ by $k_1-1$ in the last formula we obtain
\begin{multline}
\label{I2}
 I_2=\sum_{k=1}^{\left[\frac{m}{2}\right]+1}\sum_{\substack{(k_1,\cdots,k_n)\in\NN^n_0\\k_1+\dots+k_n=k, \ k_1\geq 1}}\frac{c_{m+1-k}}{2^{k-1}(k_1-1)!k_2!\dots k_n!}|x|^{2-n-2m+2(k-1)}\frac{1}{\alpha_1+1}D^{2k}_{x_1^{2k_1},\dots,x_n^{2k_n}}x_1x^\alpha\\
=\sum_{k=1}^{\left[\frac{m}{2}\right]}\sum_{\substack{(k_1,\cdots,k_n)\in\NN^n_0\\k_1+\dots+k_n=k}}\frac{c_{m+1-k}}{2^{k}k_1!k_2!\dots k_n!}|x|^{2-n-2(m+1)+2k}\frac{2k_1}{\alpha_1+1}D^{2k}_{x_1^{2k_1},\dots,x_n^{2k_n}}x_1x^\alpha +R_m,
\end{multline}
where
$$R_m=\sum_{\substack{(k_1,\cdots,k_n)\in\NN^n_0\\k_1+\dots+k_n=\left[\frac{m}{2}\right]+1}}\frac{c_{m-\left[\frac{m}{2}\right]}}{2^{\left[\frac{m}{2}\right]+1}k_1!k_2!\dots k_n!}|x|^{2-n-2m+2\left[\frac{m}{2}\right]}\frac{2k_1}{\alpha_1+1}D^{2(\left[\frac{m}{2}\right]+1)}_{x_1^{2k_1},\dots,x_n^{2k_n}}x_1x^\alpha.$$
From (\ref{eI1}) and (\ref{I2}) we get
\begin{eqnarray*}
I_1+I_2=\sum_{k=1}^{\left[\frac{m}{2}\right]}\sum_{\substack{(k_1,\cdots,k_n)\in\NN^n_0\\k_1+\dots+k_n=k}}\frac{c_{m+1-k}}{2^{k}k_1!k_2!\dots k_n!}|x|^{2-n-2(m+1)+2k}D^{2k}_{x_1^{2k_1},\dots,x_n^{2k_n}}x_1x^\alpha +R_m.
\end{eqnarray*}
It is obvious that if $m$ is an even number, then $R_m=0$. Let $m$ be an odd number, then $\left[\frac{m}{2}\right]+1=\left[\frac{m+1}{2}\right]$. Let us note that the elements of the sum $R_m$ do not vanish if and only if $(k_1,k_2,...,k_n)=(\frac{\alpha_1+1 }{2},\frac{\alpha_2}{2},\dots,\frac{\alpha_n }{2})$ and $\alpha_1+1,...,\alpha_n$ are even. Therefore, without loss of generality we may assume that $\alpha_1+1,...,\alpha_n$ are even. Then
\begin{eqnarray*}
R_m&=&\frac{c_{m-\left[\frac{m}{2}\right]}}{2^{\left[\frac{m}{2}\right]+1}(\frac{\alpha_1+1 }{2})!(\frac{\alpha_2}{2})!\dots (\frac{\alpha_n}{2})!}|x|^{2-n-2m+2\left[\frac{m}{2}\right]}\frac{2\cdot\frac{\alpha_1+1 }{2}}{\alpha_1+1}D^{2(\left[\frac{m}{2}\right]+1)}_{x_1^{\alpha_1+1},\dots,x_n^{\alpha_n}}x_1x^\alpha\\
&=&\frac{c_{m+1-\left[\frac{m+1}{2}\right]}}{2^{\left[\frac{m+1}{2}\right]}(\frac{\alpha_1+1 }{2})!(\frac{\alpha_2}{2})!\dots (\frac{\alpha_n}{2})!}|x|^{2-n-2(m+1)+2\left[\frac{m+1}{2}\right]}D^{2\left[\frac{m+1}{2}\right]}_{x_1^{\alpha_1+1},\dots,x_n^{\alpha_n}}x_1x^\alpha.
\end{eqnarray*}
Finally we can write:
\begin{multline*}
D_{x_1}D^\alpha|x|^{2-n}=c_{m}(2-n-2m)|x|^{2-n-2m-2}x_1x^\alpha\\+\sum_{k=1}^{\left[\frac{m}{2}\right]}\sum_{\substack{(k_1,\cdots,k_n)\in\NN^n_0\\k_1+\dots+k_n=k}}\frac{c_{m+1-k}}{2^{k}k_1!k_2!\dots k_n!}|x|^{2-n-2(m+1)+2k}D^{2k}_{x_1^{2k_1},\dots,x_n^{2k_n}}x_1x^\alpha +R_m\\
=\sum_{k=0}^{\left[\frac{m+1}{2}\right]}\sum_{\substack{(k_1,\cdots,k_n)\in\NN^n_0\\k_1+\dots+k_n=k}}\frac{c_{m+1-k}}{2^{k}k_1!k_2!\dots k_n!}|x|^{2-n-2(m+1)+2k}D^{2k}_{x_1^{2k_1},\dots,x_n^{2k_n}}x_1x^\alpha.
\end{multline*}
Hence
\begin{multline*}
K[D_{x_1}D^\alpha|x|^{2-n}]=\sum_{k=0}^{\left[\frac{m+1}{2}\right]}\sum_{\substack{(k_1,\cdots,k_n)\in\NN^n_0\\k_1+\dots+k_n=k}}\frac{c_{m+1-k}}{2^{k}k_1!k_2!\dots k_n!}\\
\times K[|x|^{2-n-2(m+1)+2k}D^{2k}_{x_1^{2k_1},\dots,x_n^{2k_n}}x_1x^\alpha]\\
=\sum_{k=0}^{\left[\frac{m+1}{2}\right]}\sum_{\substack{(k_1,\cdots,k_n)\in\NN^n_0\\k_1+\dots+k_n=k}}\frac{c_{m+1-k}}{2^{k}k_1!k_2!\dots k_n!}|x|^{2k}D^{2k}_{x_1^{2k_1},\dots,x_n^{2k_n}}x_1x^\alpha,
\end{multline*}
because $D^{2k}_{x_1^{2k_1},\dots,x_n^{2k_n}}x_1x^\alpha\in\mathcal{P}_{m+1-2k}(\RR^n)$ for $k=0,1,\dots,\left[\frac{m+1}{2}\right] $ and for every $(k_1,...,k_n)\in\NN_0^n$ such that $k_1+\dots+k_n=k$.   One can replace the differentiation with respect to $x_1$ by the differentiation with respect to any variable $x_j$, $j=1,\dots,n$:
\begin{multline*}
K[D_{x_j}D^\alpha|x|^{2-n}]=\sum_{k=0}^{\left[\frac{m+1}{2}\right]}\sum_{\substack{(k_1,\cdots,k_n)\in\NN^n_0\\k_1+\dots+k_n=k}}\frac{c_{m+1-k}}{2^{k}k_1!k_2!\dots k_n!}|x|^{2k}D^{2k}_{x_1^{2k_1},\dots,x_n^{2k_n}}x_jx^\alpha.
\end{multline*}
Since the expression $x_jx^\alpha$ determines any monomial of degree $m+1$, the proof of (i) is finished by the mathematical induction.

A proof of (ii) is almost the same. Since $K[D_{x_i}\log|x|]=x_i$,
we notice that $\tilde{c}_1=1$.
Here we have $n=2$, hence $\tilde{c}_{m+1}=-2m\tilde{c}_m$. Therefore we conclude that 
$$\tilde{c}_{m+1}=(-2)^mm!.$$
\end{proof}
From the proved proposition we get a following corollary, which extends the result \cite[Corollary 5.20]{A-B-R} formulated only for homogeneous harmonic polynomials to any homogeneous polynomials.
\begin{Wn}
\label{cor:1}
Let $u\in\mathcal{P}_m(\RR^n)$ with $m\geq 1$.
\begin{enumerate}
\item[(i)] If $n>2$, then
$$u(x)=\sum_{k=0}^{\left[\frac{m}{2}\right]}\frac{|x|^{2k}}{c_{m-2k}}K[u_k(D)|x|^{2-n}].$$
\item[(ii)] If $n=2$, then 
$$u(x)=\sum_{k=0}^{\left[\frac{m}{2}\right]}\frac{|x|^{2k}}{\tilde{c}_{m-2k}}K[u_k(D)\log|x|].$$
\end{enumerate}
In both cases $u_0=u$ and 
$$
u_k=\sum_{\substack{ (k_1,\dots,k_n)\in\NN_0^n\\k_1+\dots +k_n=k}}d_{m,k,k_1,\dots,k_n}D^{2k}_{x_{1}^{2k_1}\dots x_{n}^{2k_n}}u\quad\textrm{for}\quad k=1,\dots,[m/2],
$$
where $d_{m,k,k_1,\dots,k_n}$ are some constants being independent of $u$.
\end{Wn}

\begin{proof}
We shall prove (i) of the statement, a proof of (ii) is similar. By Proposition \ref{prop:1} we know that
\begin{multline}
\label{e8}
\frac{1}{c_{m}} K[u(D)|x|^{2-n}]=u+|x|^2u_1\\+\sum_{k=2}^{\left[\frac{m}{2}\right]}\sum_{\substack{(k_1,\cdots,k_n)\in\NN^n_0\\k_1+\dots+k_n=k}}\frac{c_{m-k}}{2^{k}k_1!\dots k_n!c_{m}}|x|^{2k}D^{2k}_{x_1^{2k_1},\dots,x_n^{2k_n}}u,
\end{multline}
where
$$u_1=\sum_{i=1}^n\frac{c_{m-1}}{2c_{m}} D^2_{x_i^2}u\in\mathcal{P}_{m-2}(\RR^n).$$

Using again Proposition \ref{prop:1} to the function $u_1$, we get:
\begin{eqnarray*}
&&K[u_1(D)|x|^{2-n}]\\
&=& c_{m-2}u_1+\sum_{k=1}^{\left[\frac{m-2}{2}\right]}\sum_{k_1+...+k_n=k}\frac{c_{m-2-k}}{2^kk_1!...k_n!}|x|^{2k}D^{2k}_{x_1^{2k_1}x_2^{2k_2}...x_n^{2k_n}}u_1\\
&=& c_{m-2}u_1+\sum_{k=1}^{\left[\frac{m-2}{2}\right]}\sum_{k_1+...+k_n=k}\frac{c_{m-2-k}c_{m-1}}{2^{k+1}k_1!...k_n!c_{m}}|x|^{2k}D^{2k}_{x_1^{2k_1}x_2^{2k_2}...x_n^{2k_n}}\left(\sum_{i=1}^n D^2_{x_i^2}u\right)\\
&=& c_{m-2}u_1+\sum_{k=1}^{\left[\frac{m-2}{2}\right]}\sum_{k_1+...+k_n=k+1}\frac{c_{m-2-k}c_{m-1}(k+1)}{2^{k+1}k_1!...k_n!c_{m}}|x|^{2k}D^{2(k+1)}_{x_1^{2k_1}x_2^{2k_2}...x_n^{2k_n}}u\\
&=& c_{m-2}u_1+\sum_{k=2}^{\left[\frac{m-2}{2}\right]+1}\sum_{k_1+...+k_n=k}\frac{c_{m-1-k}c_{m-1}k}{2^{k}k_1!...k_n!c_{m}}|x|^{2(k-1)}D^{2k}_{x_1^{2k_1}x_2^{2k_2}...x_n^{2k_n}}u\\
&=& c_{m-2}u_1+\sum_{k=2}^{\left[\frac{m}{2}\right]}\sum_{k_1+...+k_n=k}\frac{c_{m-1-k}c_{m-1}k}{2^{k}k_1!...k_n!c_{m}}|x|^{2(k-1)}D^{2k}_{x_1^{2k_1}x_2^{2k_2}...x_n^{2k_n}}u.
\end{eqnarray*}
So
\begin{equation*}
u_1=\frac{1}{c_{m-2}}K[u_1(D)|x|^{2-n}]-\sum_{k=2}^{\left[\frac{m}{2}\right]}\sum_{k_1+...+k_n=k}\frac{c_{m-1-k}c_{m-1}k}{2^{k}k_1!...k_n!c_{m}c_{m-2}}|x|^{2(k-1)}D^{2k}_{x_1^{2k_1}x_2^{2k_2}...x_n^{2k_n}}u.
\end{equation*}
Putting the above to (\ref{e8})  we obtain 
\begin{multline*}
\frac{1}{c_{m}}K[u_0(D)|x|^{2-n}]-\frac{|x|^2}{c_{m-2}} K[u_1(D)|x|^{2-n}]\\
=u+\sum_{k=2}^{\left[\frac{m}{2}\right]}\sum_{k_1+...+k_n=k}\tilde{d}_{m,k,k_1,...,k_n}|x|^{2k}D^{2k}_{x_1^{2k_1}x_2^{2k_2}...x_n^{2k_n}}u,
\end{multline*}
where $$\tilde{d}_{m,k,k_1,...,k_n}=\frac{c_{m-k}c_{m-2}-c_{m-1-k}c_{m-1}k}{2^{k}k_1!...k_n!c_{m}c_{m-2}}.$$
By the linearity of the Kelvin transform we may write after changing appropriate constants
\begin{multline*}
\frac{1}{c_{m}}K[u_0(D)|x|^{2-n}]+\frac{|x|^2}{c_{m-2}} K[u_1(D)|x|^{2-n}]\\
=u+\sum_{k=2}^{\left[\frac{m}{2}\right]}\sum_{k_1+...+k_n=k}\tilde{d}_{m,k,k_1,...,k_n}|x|^{2k}D^{2k}_{x_1^{2k_1}x_2^{2k_2}...x_n^{2k_n}}u.
\end{multline*}
From the last equation we have
\begin{multline*}
\frac{1}{c_{m}}K[u_0(D)|x|^{2-n}]+\frac{|x|^2}{c_{m-2}} K[u_1(D)|x|^{2-n}]\\
=u+|x|^4u_2+\sum_{k=3}^{\left[\frac{m}{2}\right]}\sum_{k_1+...+k_n=k}\tilde{d}_{m,k,k_1,...,k_n}|x|^{2k}D^{2k}_{x_1^{2k_1}x_2^{2k_2}...x_n^{2k_n}}u,
\end{multline*}
where
$$u_2=\sum_{k_1+...+k_n=2}\frac{c^2_{m-2}-2c_{m-3}c_{m-1}}{4k_1!...k_n!c_{m}c_{m-2}}D^{4}_{x_1^{2k_1}x_2^{2k_2}...x_n^{2k_n}}u\in\mathcal{P}_{m-4}(\RR^n).$$
Continuing the process  with respect to the function $u_2$ and then with respect to the functions $u_3,\dots,u_{[m/2]}$ we conclude finally that
$$u(x)=\sum_{k=0}^{\left[\frac{m}{2}\right]}\frac{|x|^{2k}}{c_{m-2k}}K[u_k(D)|x|^{2-n}].$$
\end{proof}

\begin{Wn}
\label{cor:2}
Let $m\geq 2p$ and $u\in\mathcal{P}_m(\RR^n)$. 
\begin{enumerate}
\item[(i)] If $n>2$, then
$$\sum_{k=0}^{p-1}\frac{|x|^{2k}}{c_{m-2k}}K[u_k(D)|x|^{2-n}]\in\mathcal{H}_m^p(\RR^n)$$
and
$$\sum_{k=0}^{p-1}\frac{|x|^{2k}}{c_{m-2k}}K[u_k(D)|x|^{2-n}]=u-|x|^{2p}v.$$
\item[(ii)] If $n=2$, then
$$\sum_{k=0}^{p-1}\frac{|x|^{2k}}{\tilde{c}_{m-2k}}K[u_k(D)\log |x|]\in\mathcal{H}_m^p(\RR^n)$$
and
$$\sum_{k=0}^{p-1}\frac{|x|^{2k}}{\tilde{c}_{m-2k}}K[u_k(D)\log|x|]=u-|x|^{2p}v$$
\end{enumerate}
where in the both cases $v\in\mathcal{P}_{m-2p}(\RR^n)$ and
\begin{equation}
\label{eq:v1}
v=\sum_{k=p}^{\left[\frac{m}{2}\right]}\sum_{\substack{(k_1,...,k_n)\in\NN_0\\k_1+...+k_n=k}}d_{m,k,k_1,...,k_n}|x|^{2(k-p)}D^{2k}_{x_1^{2k_1},\dots,x_{n}^{2k_n}} u,
\end{equation}
with some constants $d_{m,k,k_1,...,k_n}\in\CC$.
\end{Wn}
\begin{proof}
We will prove (i). Repeating the proof of Corollary \ref{cor:1} for $k=0,\dots,p-1$
we see that
$$
u(x)=\sum_{k=0}^{p-1}\frac{|x|^{2k}}{c_{m-2k}}K[u_k(D)|x|^{2-n}]+|x|^{2p}v(x),
$$
where $v(x)$ is given by (\ref{eq:v1}). 
Moreover, since $u_k\in\mathcal{P}_{m-2k}(\RR^n)$, by \cite[Lemma 5.15]{A-B-R} we get
$K[u_k(D)|x|^{2-n}]\in\mathcal{H}_{m-2k}(\RR^n)$. Hence by Lemma \ref{Almansi} we conclude that
$$\sum_{k=0}^{p-1}\frac{|x|^{2k}}{c_{m-2k}}K[u_k(D)|x|^{2-n}]\in\mathcal{H}_m^p(\RR^n).$$
\end{proof}

Now we are ready to state the main result about the orthogonal projection of homogeneous polynomials onto the space of homogeneous polyharmonic polynomials. The presented result is a polyharmonic version of \cite[Theorem 5.18]{A-B-R}.
\begin{Tw}
\label{th:1}
Let $m\geq 2p$ and $u\in\mathcal{P}_m(\RR^n)$. Then the mapping $Q\colon\mathcal{P}_m(\RR^n)\rightarrow\mathcal{H}_m^p(\RR^n)$ of the form
$$
Q[u](x)=\systeme*{\sum_{k=0}^{p-1}\frac{|x|^{2k}}{c_{m-2k}}K[u_k(D)|x|^{2-n}] \  \ \ \ for \ n>2{,},\sum_{k=0}^{p-1}\frac{|x|^{2k}}{\tilde{c}_{m-2k}}K[u_k(D)\log |x|] \  \ \ \ for  \ n=2}
$$
is the canonical projection of $\mathcal{P}_m(\RR^n)$ onto $\mathcal{H}_m^p(\RR^n)$, which extends in a natural way to the canonical projection $\tilde{Q}$ of $\mathcal{P}_m(\CC^n)$ onto $\mathcal{H}_m^p(\CC^n)$.
\end{Tw}
\begin{proof}
We shall prove the case $n>2$. A proof of the case $n=2$ is similar.
By Corollary \ref{cor:2}, we have a unique decomposition of $u\in\mathcal{P}_m(\RR^n)$ given by
\begin{equation*}
u(x)=\sum_{k=0}^{p-1}\frac{|x|^{2k}}{c_{m-2k}}K[u_k(D)|x|^{2-n}]+|x|^{2p}v(x)
\end{equation*}
for $v(x)$ satisfying (\ref{eq:v1}), where $\sum_{k=0}^{p-1}\frac{|x|^{2k}}{c_{m-2k}}K[u_k(D)|x|^{2-n}]\in\mathcal{H}_m^p(\RR^n)$ and the last term $|x|^{2p}v(x)$ does not belong to $\mathcal{H}_m^p(\RR^n)$. Hence the mapping $Q$ is the canonical projection of $\mathcal{P}_m(\RR^n)$ onto $\mathcal{H}_m^p(\RR^n)$.

Since every polynomial $u(x)$ for $x\in\RR^n$ extends in a natural way to its complexification $u(z)$ for $z\in\CC^n$,
the mapping $Q$ extends in a natural way to the canonical projection $\tilde{Q}$ of $\mathcal{P}_m(\CC^n)$ onto $\mathcal{H}_m^p(\CC^n)$.
\end{proof}

\section{Projection of the complex homogeneous polynomials onto the spherical polyharmonics}
In this section we derive the projection of the complex homogeneous polynomials onto the space of spherical polyharmonics. A spherical polyharmonic is a restriction of the homogeneous polynomial to the union of rotated spheres $\widehat{S}_p:=\bigcup_{l=0}^{p-1}e^{\frac{l\pi i}{p}}S$. The motivation to study polyharmonic functions on the union of rotated balls $\widehat{B}_p:=\bigcup_{l=0}^{p-1}e^{\frac{l\pi i}{p}}B$ comes from the Pizzetti-type formula for the operator $\Delta^p$ given
in \cite{L0} and \cite{S-M} (see also \cite{G-M} and \cite{G-M-2}). By this formula, the integral mean of $u$ over the rotated spheres $x+\bigcup_{l=0}^{p-1}e^\frac{l\pi i}{p}S(0,r)$ given by
$$M_{\Delta^p}(u;x,r):=\frac{1}{p\omega_n}\sum_{l=0}^{p-1}\int\limits_{S}u(x+e^{\frac{l\pi i}{p}}r\zeta)\,dS(\zeta)
$$
has the expansion
$$
 M_{\Delta^p}(u;x,r)
 =\sum_{j=0}^{\infty}\frac{\Delta^{pj}u(x)}{4^{pj}(n/2)_{pj}(pj)!}r^{2pj},
 $$
where $(a)_k:=a(a+1)\cdots(a+k-1)$, for $k\in\NN$, is the Pochhammer symbol and $\omega_n $ denotes the area of the unit sphere in $ {\RR}^n $.
Hence, the mean value property holds  for every
 polyharmonic function $u$ of order $p$ on the closed rotated balls $x_0+\bigcup_{k=0}^{p-1}e^\frac{k\pi i}{p}\overline{B}(0,r)$:
 $$u(x)=\frac{1}{p \omega_n}\sum_{l=0}^{p-1}\int\limits_{S}u(x+e^\frac{l\pi i}{p}r\zeta)\,dS(\zeta).$$
In particular, it means that the value $u(0)$ of the polyharmonic function $u$ is uniquely determined by the boundary values on the rotated unit
spheres $\widehat{S}_p$:
$$u(0)=\frac{1}{p \omega_n}\sum_{k=0}^{p-1}\int\limits_{S}u(e^\frac{k\pi i}{p}r\zeta)\,dS(\zeta).$$

 \begin{Df}[{\cite[Definition 1]{G-M-2}}]
\label{D1}
The restriction to the set $\widehat{S}_p$ of an element of $\mathcal{H}_m^p(\CC^n)$
is called a \emph{spherical polyharmonic of degree $m$ and order $p$}.

The set of spherical polyharmonics is denoted by $\mathcal{H}_m^p( \widehat{S}_p ) $.
\end{Df}

The spherical polyharmonics of order $1$ are called \emph{spherical harmonics} and their space is denoted by $\mathcal{H}_m(S):=\mathcal{H}_m^1(S)$
(see \cite[Chapter 5]{A-B-R}). Analogously we write $\mathcal{H}_{m}(\CC^n)$ instead of $\mathcal{H}_{m}^1(\CC^n)$.

Let us consider the Hilbert space $ L^2 ( \widehat{S}_p )  $ of square-integrable functions on $\widehat{S}_p$ with the inner product defined by
\begin{equation}
\label{inner product 1}
\left\langle f,g\right\rangle _{ \widehat{S}_p }:=\frac{1}{p} \int\limits_S \sum_{j=0}^{p-1}f(e^{\frac{j\pi i}{p}} \zeta )
\overline{g(e^{\frac{j\pi i}{p}} \zeta )}\,d\sigma (\zeta),
\end{equation}
where $d\sigma$ is a normalized surface-area measure on the unit sphere $S$.

It is known that  $\mathcal{H}^p_m(\widehat{S}_p)$ is finite dimensional
{\cite[Proposition 4]{G-M-2}} and $h_m^p:=\dim\mathcal{H}^p_m(\CC^n)=\mathcal{H}^p_m(\widehat{S}_p)$, moreover

\begin{eqnarray*}
 h_m^p&=&\left\{
  \begin{array}{lll}
    \binom{n+m-1}{n-1} & \textrm{for} & m<2p,\\
    \binom{n+m-1}{n-1}-\binom{n+m-2p-1}{n-1} & \textrm{for} & m\geq 2p.
  \end{array}
  \right.
\end{eqnarray*}

\begin{Le}[{\cite[Theorem 1]{G-M-2}}]
\label{L3}
The space $L^2 ( \widehat{S}_p )$ is the \emph{direct sum} of spaces $ \mathcal{H}^p_m(\widehat{S}_p)$, which is written as
\begin{equation*}
L^2 ( \widehat{S}_p )=\bigoplus_{m=0}^{\infty} \mathcal{H}_m^p ( \widehat{S}_p ).
\end{equation*}
It means that
\begin{enumerate}
\item[(i)] $\mathcal{H}^p_m(\widehat{S}_p)$ is a closed subspace of $L^2 ( \widehat{S}_p )$ for every $m$.
\item[(ii)] $\mathcal{H}^p_m(\widehat{S}_p)$ is orthogonal to $\mathcal{H}^p_k(\widehat{S}_p)$ if $m\neq k$.
\item[(iii)] For every $ x \in L^2 ( \widehat{S}_p )$ there exist $ x_m \in \mathcal{H}^p_m(\widehat{S}_p)$, $m=0,1,...$,  such that $ x=x_0+x_1+x_2+\dots$, where the sum is convergent in the norm of $L^2 ( \widehat{S}_p )$.
\end{enumerate}
\end{Le}

Now we are ready to formulate the counterpart of Theorem \ref{th:1} for spherical polyharmonics.
\begin{Tw}
Let $m\geq 2p$ and $u\in\mathcal{P}_m(\CC^n)$. Then the mapping $\tilde{Q}|_{\widehat{S}_p}\colon\mathcal{P}_m(\CC^n)|_{\widehat{S}_p}\rightarrow\mathcal{H}_m^p(\widehat{S}_p)$ is the orthogonal projection of $\mathcal{P}_m(\CC^n)|_{\widehat{S}_p}$ onto $\mathcal{H}_m^p(\widehat{S}_p)$ given by
$$
\tilde{Q}|_{\widehat{S}_p}[u|_{\widehat{S}_p}](z)=\systeme*{e^{iml\pi/p}\sum_{k=0}^{p-1} u_k(D)|x|^{2-n}/c_{m-2k} \  \ \ \ \textrm{for} \ n>2 , e^{iml\pi/p}\sum_{k=0}^{p-1}u_k(D)\log |x|/\tilde{c}_{m-2k} \  \ \ \ \textrm{for}  \ n=2},
$$
where $z\in\widehat{S}_p\subset\CC^n$ and $z=e^{il\pi/p}x$ for some $l\in\{0,\dots,p-1\}$ and $x\in S$.
\end{Tw}
\begin{proof}
We shall prove the theorem for $n>2$. Again, by Corollary \ref{cor:2}, we have a unique decomposition of $u\in\mathcal{P}_m(\RR^n)$ given by
\begin{equation}
\label{eq:decomposition}
u(x)=\sum_{k=0}^{p-1}\frac{|x|^{2k}}{c_{m-2k}}K[u_k(D)|x|^{2-n}]+|x|^{2p}v(x)
\end{equation}
for $v(x)$ satisfying (\ref{eq:v1}), where $\sum_{k=0}^{p-1}\frac{|x|^{2k}}{c_{m-2k}}K[u_k(D)|x|^{2-n}]\in\mathcal{H}_m^p(\RR^n)$ and the last term $|x|^{2p}v(x)$ does not belong to $\mathcal{H}_m^p(\RR^n)$. We restrict both sides of the complexification of (\ref{eq:decomposition}) to $e^{il\pi/p}S\subset\CC^n$ for some $l\in\{0,\dots,p-1\}$. Since the Kelvin transform preserves points of the unit sphere and the elements of the sum (\ref{eq:decomposition}) are homogeneous of degree $m$, we get
$$
u(e^{il\pi/p}x)=e^{iml\pi/p}\sum_{k=0}^{p-1} u_k(D)|x|^{2-n}/c_{m-2k}+|x|^{2p}v(e^{il\pi/p}x).
$$
By \cite[Propositions 2 and 5]{G-M-2} $u|_{\widehat{S}_p}$ is orthogonal to $\mathcal{H}_m^p(\widehat{S}_p)$ in $L^2(\widehat{S}_p)$.
\end{proof}

\section{Basis of the space $\mathcal{H}^p_m(\CC^n)$}
In this section we want to provide a basis for the space $\mathcal{H}^p_m(\CC^n)$, the
space of homogeneous polyharmonic  polynomials of order $p$ and  degree $m$. For the case $p=1$ and
$n=2$ it is known that $\mathcal{H}_m(\RR^2)$ is the complex linear span of $\{z^m,\overline{z^m}\}$, and for $n>2$ the basis of $\mathcal{H}_m(\RR^n)$ is given in \cite[Theorem 5.25]{A-B-R}. We now prove the following result:
\begin{Tw}
\label{th:basis}
Set $v_{1,k}=\RE(x_1+ix_2)^{m-2k}|x|^{2k}, \ v_{2,k}=\IM(x_1+ix_2)^{m-2k}|x|^{2k}$. Then 
$$\Big\{v_{1,k},v_{2,k}\colon k=0,\dots,p-1\Big\}
$$ is a vector space basis of $\mathcal{H}_m^p(\RR^2)$.

Set $v_{k,\alpha(k)}=|x|^{2k}K[D^{\alpha(k)}|x|^{2-n}]$. Then the set
\begin{multline*}
\Big\{v_{k,\alpha(k)}\colon \alpha(k)\in\NN^n_0,\, |\alpha(k)|=m-2k,\, \alpha(k)_1\leq 1, \  k=0,\dots,p-1\Big\}
\end{multline*}
is a vector space basis of $\mathcal{H}_m^p(\RR^n)$ with $n>2$.
\end{Tw}
 \begin{proof}
Let $n=2$. It is known that (see \cite[page 82]{A-B-R})
$$\mathcal{H}_m(\RR^2)=\lin \{v_{1,0}, v_{2,0}\},$$
and (see \cite[Theorem 5.25]{A-B-R})
$$\mathcal{H}_m(\RR^n)=\lin \{K[D^{\alpha}|x|^{2-n}]:|\alpha|{=m} {,} \ \alpha_1\leq 1\} \ \textrm{for} \ n>2.$$

So, by Lemma \ref{Almansi} we get that every polynomial $q\in \mathcal{H}^p_m(\CC^n)$ can be written as a linear combination of the vectors $v_{1,k}, v_{2,k}$ and $v_{k,\alpha(k)}$ when $n=2$ and $n>2$, respectively. Moreover, by Lemma \ref{Almansi} the number of these vectors is equal to dimension $h^p_m$ (see also \cite[Proposition 9]{G-M-2}). Hence these vectors create a basis of $\mathcal{H}_m^p(\RR^2)$ and $\mathcal{H}_m^p(\RR^n)$ ($n>2$), respectively.
 \end{proof}
Consequently, we conclude that
\begin{Wn}
If $n=2$, then the set
\begin{multline*}
\Big\{w_{1,k}(z)=e^{iml\pi/p}\cos(m-2k)\theta,\ w_{2,k}(z)=e^{iml\pi/p}\sin(m-2k)\theta\colon\\
k=0,\dots,p-1,\ \textrm{for}\ z=e^{il\pi/p}(\cos\theta,\sin\theta),\ l\in\{0,\dots,p-1\},\ \theta\in[0,2\pi)\Big\}
\end{multline*}
is a vector space basis of $\mathcal{H}_m^p(\widehat{S}_p)$. 

Moreover for $n>2$ a vector space basis of $\mathcal{H}_m^p(\widehat{S}_p)$ is given by the set
\begin{multline*}
 \Big\{w_{k,\alpha(k)}(z)=e^{iml\pi/p}D^{\alpha(k)}|x|^{2-n}\colon 
 |\alpha(k)|=m-2k,\ \alpha(k)_1\leq 1,\\ k=0,\dots,p-1,\ \textrm{for}\ z=e^{il\pi/p}x,\ l\in\{0,\dots,p-1\},\ x\in S\Big\}.
\end{multline*}
\end{Wn}
\begin{proof}
 By Theorem \ref{th:basis} a basis of $\mathcal{H}_m^p(\RR^2)$ is given by following polynomials of two real variables $x=(x_1,x_2)\in\RR^2$ 
 \begin{equation*}
  v_{1,k}(x)=\RE(x_1+ix_2)^{m-2k}|x|^{2k}\quad\textrm{and}\quad v_{2,k}(x)=\IM(x_1+ix_2)^{m-2k}|x|^{2k},
 \end{equation*}
where $k=0,\dots,p-1$.

Using the complexification of polynomials $v_{j,k}(x)$ ($j=1,2$, $k=0,\dots,l-1$) we conclude that the set
$$\{v_{1,k}(z), v_{2,k}(z)\colon k=0,\dots,p-1,\, z=(z_1,z_2)\in\CC^2\}$$
is a vector space basis of $\mathcal{H}_m^p(\CC^2)$.

It means that a vector space basis of $\mathcal{H}_m^p(\widehat{S}_p)$ is
given by 
$$\{v_{1,k}|_{\widehat{S}_p},v_{2,k}|_{\widehat{S}_p}\colon k=0,\dots,p-1\}.$$

Let $w_{1,k}=v_{1,k}|_{\widehat{S}_p}, \ w_{2,k}=v_{2,k}|_{\widehat{S}_p}.$ Take any $z\in\widehat{S}_p\subseteq\CC^2$ and observe that $z=e^{il\pi/p}(\cos \theta, \sin\theta)$ for some $l\in\{0,\dots,p-1\}$ and $\theta\in[0,2\pi)$.
For such $z$ we get
$$w_{1,k}(z)=e^{ilm\pi/p}w_{1k}(\cos\theta,\sin\theta)=e^{ilm\pi/p}\RE e^{i(m-2k)\theta}=e^{ilm\pi/p}\cos(m-2k)\theta$$
and analogously
$w_{2,k}(z)=e^{ilm\pi/p}\sin(m-2k)\theta$.

Analogously, for $n>2$ we define polynomials in $\mathcal{H}_m^p(\RR^n)$
$$
v_{k,\alpha(k)}(x)=|x|^{2k}K[D^{\alpha(k)}|x|^{2-n}],
$$
where $|\alpha(k)|=m-2k$, $\alpha(k)_1\leq 1$ and $k=0,\dots,p-1$. Using their complexification $v_{k,\alpha(k)}(z)$ for $z\in\CC^n$, by Theorem \ref{th:basis} we conclude that the set
$$\{v_{k,\alpha(k)}(z)\colon |\alpha(k)|=m-2k,\ \alpha(k)_1\leq 1,\ k=0,\dots,p-1\}$$
is a vector space basis of $\mathcal{H}_m^p(\CC^n)$.

Hence a vector space basis of $\mathcal{H}_m^p(\widehat{S}_p)$ is given by
$$\{v_{k,\alpha(k)}|_{\widehat{S}_k}\colon |\alpha(k)|=m-2k,\ \alpha(k)_1\leq 1,\ k=0,\dots,p-1\}.$$
Any $z\in \widehat{S}_p\subset \CC^n$ may be written as
$z=e^{il\pi/p}x$ for some $l\in\{0,\dots,p-1\}$ and $x\in S$. Let $w_{k,\alpha(k)}=v_{k,\alpha(k)}|_{\widehat{S}_k}.$ Then
$$w_{k,\alpha(k)}(z)=e^{ilm\pi/p}v_{k,\alpha(k)}(x)=e^{ilm\pi/p}D^{\alpha(k)}|x|^{2-n}.$$

\end{proof}

\section{Convergence of the orthogonal series}
We shall prove (see Theorem 3.1 in \cite{F-V}):
\begin{Tw}
Let $\{e_{m,1},\dots,e_{m,h^p_m}\}$ be an orthonormal basis of $\mathcal{H}_m^p(\CC^n)$ with respect to the inner product (\ref{inner product 1}). Then the series 
$$f(x)=\sum_{m=0}^\infty\sum_{j=1}^{h_m^p}a_{m,j}e_{m,j}(x)$$
converges absolutely and uniformly on compact subsets of $\widehat{B}_p(R)$, where
\begin{equation}
\label{eq:R}
 R^{-1}=\limsup_{m\rightarrow\infty}\left(||a_{m}||\right)^{1/m} \  {and} \ ||a_m||^2=\sum_{j=1}^{h^p_m}|a_{m,j}|_{\CC}^2.
\end{equation}
\end{Tw}

\begin{proof}
Let $x\in\widehat{B}_p(r)$ for some $r>0$. We have $e_{m,j}(x)=r^me_{m,j}(\frac{x}{r})$, hence by the Cauchy-Schwarz inequality
$$\left|\sum_{m=0}^\infty\sum_{j=1}^{h_m^p}a_{m,j}e_{m,j}(x)\right|_{\CC}\leq \sum_{m=0}^\infty r^m \left|\sum_{j=1}^{h_m^p}a_{m,j}e_{m,j}\left(\frac{x}{r}\right)\right|_{\CC}\leq$$
$$\leq \sum_{m=0}^\infty r^m \left(\sum_{j=1}^{h_m^p}|a_{m,j}|_{\CC}^2\right)^{1/2}\left(\sum_{j=1}^{h_m^p} \left|e_{m,j}\left(\frac{x}{r}\right)\right|^2_{\CC}\right)^{1/2}\leq \sum_{m=0}^\infty r^m||a_{m}||\sqrt{h^p_m}.$$
Since $e_{m,j}(x)$ are homogeneous polynomials of order $m>0$,
by Lemma~\ref{St13} and Lemma~\ref{Lemat 5} in the last inequality we can use the estimation
$$
\sum_{j=1}^{h_m^p} \left|e_{m,j}\left(\frac{x}{r}\right)\right|^2_{\CC}\leq \sum_{j=1}^{h_m^p} \left|e_{m,j}\left(\frac{x}{|x|_{\CC}}\right)\right|^2_{\CC}=Z^p_m(\frac{x}{|x|_{\CC}},\frac{x}{|x|_{\CC}})=h^p_m.
$$
By the Cauchy-Hadamard theorem we get the  formula for the radius~$R$:
$$R^{-1}=\limsup_{m\rightarrow\infty}\left(\sqrt{h^p_m}||a_{m}||\right)^{1/m}.$$
 Moreover
$$h^p_m= {{n+m-1}\choose{n-1}}-{{n+m-2p-1}\choose{n-1}},$$
so
$$1\leq h^p_m\leq {{n+m-1}\choose{n-1}}\leq (n+m)^n,$$
hence 
$$\lim_{m\rightarrow\infty}\left(\sqrt{h^p_m}\right)^{\frac{1}{m}}=1,$$
which completes the proof.
\end{proof}
\begin{Uw}
 The same result holds, if we replace in (\ref{eq:R}) $\ell_2$-norm $||a_m||$ by an equivalent norm, such as $\ell_1$-norm $||a_m||_1=\sum_{j=1}^{h^p_m}|a_{m,j}|_{\CC}$ or $\ell_{\infty}$-norm $||a_m||_{\infty}=\max\{|a_{m,j}|_{\CC}\colon j=1,\dots,h^p_m\}$.
\end{Uw}

\begin{Uw}
A similar result given in \cite{F-V} is improved in \cite{R} by the explicit formula for the zonal polyharmonics studied in the mentioned paper. 
\end{Uw}

\section{Relation with the Fischer inner product}
This section we start with zonal polyharmonics, let us recall some needed information about them. 

By  Lemma \ref{L3}  we may treat $\mathcal{H}_m^p (\widehat{S}_p)$ as a Hilbert space  with the inner product (\ref{inner product 1}) induced
from $L^2( \widehat{S}_p )$. Let $ \eta \in \widehat{S}_p$ be a fixed point, by classical argument we conclude that there exists a unique $Z_m^p(\cdot,\eta)\in \mathcal{H}_m^p (\widehat{S}_p)$ such that (see {\cite[formula (9)]{G-M-2}})
\begin{equation}
\label{e'2}
q(\eta)=\left\langle q,Z^p_m(\cdot,\eta)\right\rangle_{  \widehat{S}_p }\quad\textrm{for every}\quad q\in
\mathcal{H}^p_m(\widehat{S}_p).
\end{equation}

\begin{Df} [{\cite[Definition 2]{G-M-2}}]
\label{D2}
The function $Z^p_m(\cdot,\eta)$ satisfying (\ref{e'2}) is called a \emph{zonal polyharmonic}
of degree $m$ and of order $p$ with a pole~$\eta$. 
\end{Df}
\begin{Uw}
 By the property (\ref{e'2}), the function $(\zeta,\eta)\mapsto Z^p_m(\zeta,\eta)$ defined on $\widehat{S}_p\times\widehat{S}_p$ is also called a \emph{reproducing kernel for the space $\mathcal{H}^p_m(\widehat{S}_p)$}.
\end{Uw}

Zonal polyharmonics of order $p=1$ are called \emph{zonal harmonics} and  we  denote them by
$Z_m(\cdot,\eta)$ instead of $Z^1_m(\cdot,\eta)$ for $\eta \in S$.

\begin{Le}{\cite[Proposition 7 and Theorem 2]{G-M-2}}
\label{St13}
Let $\zeta,\eta\in\widehat{S}_p$. Then 
\begin{enumerate}
\item[(i)] $Z^p_m(\zeta,\eta)=\sum_{k=1}^{h_m^p}e_k(\zeta)\overline{e_k(\eta)}$, where $\{e_1,\dots,e_{h^p_m}\}$ is an orthonormal basis of $\mathcal{H}^p_m(\widehat{S}_p)$.
\item[(ii)] $Z_m^p(\zeta,\eta)=\sum_{k=0}^{\min\{[m/2],p-1\}}|\zeta|^{2k}|\overline{\eta}|^{2k}Z_{m-2k}(\zeta,\eta)$.
\end{enumerate}
\end{Le}

\begin{Le}{\cite[Proposition 8]{G-M-2}}
\label{Lemat 5}
If $\eta\in \widehat{S}_p$, then
\begin{enumerate}
\item[(i)] $Z^p_m(\eta,\eta)=\dim\mathcal{H}^p_m(\CC^n)$.
\item[(ii)] $|Z_m^p(\zeta,\eta)|_\CC\leq \dim\mathcal{H}^p_m(\CC^n)$ for all $\zeta\in \widehat{S}_p.$
\end{enumerate}
\end{Le}

H. Render in \cite{R} studies the space of all homogeneous polynomials $\mathcal{H}_m^p(\RR^n)$ of degree $m$ and polyharmonic order $p$ endowed with 
the \emph{Fischer inner product} defined by
\begin{equation}
\label{eq:Fischer}
\langle P,Q \rangle_F:=\sum_{|\alpha|\leq N}\alpha!c_{\alpha}\overline{d_{\alpha}}
\end{equation}
for polynomials $P(x)=\sum_{|\alpha|\leq N}c_{\alpha}x^{\alpha}$ and $Q(x)=\sum_{|\alpha|\leq N} d_{\alpha} x^{\alpha}$ in $\RR^n$.

In this section we discuss the relation between the space of polyharmonic polynomials $\mathcal{H}_m^p(\CC^n)$ endowed with the Fischer inner product (\ref{eq:Fischer}), and endowed with the integral inner product (\ref{inner product 1}) inherited from the space $L^2(\widehat{S}_p)$.

In the special case $p=1$ of harmonic polynomials $\mathcal{H}_m(\CC^n)$ we have the following direct relation between these two inner products 
\begin{Le}[{\cite[Theorem 5.14]{A-B-R}, see also \cite[Theorem 2.1]{R}}]
\label{Fischer}
For $u,v\in\mathcal{H}_m(\CC^n)$ we have 
\begin{equation*}
 \left\langle u,v\right\rangle_F=C\left\langle u,v\right\rangle_{\widehat{S}_1},
\end{equation*}
where $C=C(n,m)=n(n+2)\cdots (n+2m-2)$ is a constant which depends only on the dimension $n$ of the space $\CC^n$, and on the degree $m$ of the polynomials in $\mathcal{H}_m(\CC^n)$. 

\end{Le}
It means that these both inner products are the same on the space
$\mathcal{H}_m(\CC^n)$ up to the constant $C=C(n,m)$.
We cannot extend such type result to the space
$\mathcal{H}^p_m(\CC^n)$. To this end take any $u,v\in \mathcal{H}^p_m(\CC^n)$ and observe that
by the Almansi theorem (Lemma \ref{Almansi}) we may write
\begin{equation}
\label{eq:u}
 u(x)=\sum_{j=0}^{p-1} u_j(x)|x|^{2j}\quad \textrm{with} \quad u_j(x)=\sum_{|\alpha|=m-2j}a_{j,\alpha}x^{\alpha} \in \mathcal{H}_{m-2j}(\CC^n).
\end{equation}
Analogously
\begin{equation}
\label{eq:v}
 v(x)=\sum_{j=0}^{p-1} v_j(x)|x|^{2j}\quad \textrm{with} \quad v_j(x)=\sum_{|\alpha|=m-2j}b_{j,\alpha}x^{\alpha} \in \mathcal{H}_{m-2j}(\CC^n).
\end{equation}

Calculating the Fischer and integral products we get
\begin{Tw}
 \label{Prop:difference}
 Let $u,v\in \mathcal{H}^p_m(\CC^n)$ be given by (\ref{eq:u}) and (\ref{eq:v}), respectively. Then
 \begin{equation}
 \label{eq:S}
  \langle u, v \rangle_{\widehat{S}_p}=\sum_{j=0}^{p-1}\sum_{|\alpha|=m-2j}a_{j,\alpha}\overline{b_{j,\alpha}}\frac{\alpha!}{n(n+2)\cdots (n+2m-2(2j+1))}
 \end{equation}
 and
\begin{equation}
 \label{eq:F} 
 \langle u, v \rangle_{F}=\sum_{j,k=0}^{p-1}\sum_{|\beta|=j}\sum_{|\tilde{\beta}|=k}
 \sum_{|\alpha|=m-2j \atop {|\tilde{\alpha}|=m-2k \atop \alpha+2\beta=\tilde{\alpha}+2\tilde{\beta}}}
 (\alpha+2\beta)!
 \frac{j!}{\beta!}\frac{k!}{\tilde{\beta}!}
 a_{j,\alpha}\overline{b_{k,\tilde{\alpha}}}. 
\end{equation}
\end{Tw}
\begin{proof}
To prove (\ref{eq:S}), observe that by \cite[Proposition 5]{G-M-2} $\langle u_j, v_k\rangle_{\widehat{S}_p}=0$ for $j\neq k$. Using this observation,  Lemma \ref{Fischer} and (\ref{eq:Fischer}) we conclude that
\begin{multline*}
 \langle u,v\rangle_{\widehat{S}_p}=\sum_{j,k=0}^{p-1}\langle u_j(x)|x|^{2j}, v_k(x)|x|^{2k}\rangle_{\widehat{S}_p}=
 \sum_{j=0}^{p-1}\langle u_j, v_j\rangle_{\widehat{S}_1}\\
 =\sum_{j=0}^{p-1}\sum_{|\alpha|=m-2j}a_{j,\alpha}\overline{b_{j,\alpha}}\frac{\alpha!}{n(n+2)\cdots (n+2m-2(2j+1))}.
\end{multline*}
To derive (\ref{eq:F}), first we calculate
\begin{equation*}
 |x|^{2j}=(x_1^2+\cdots+x_n^2)^j=\sum_{\beta_1+\cdots+\beta_n=j}
 \frac{j!}{\beta_1!\cdots\beta_n!}x_1^{2\beta_1}\cdots x_n^{2\beta_n}=\sum_{|\beta|=j}\frac{j!}{\beta!}x^{2\beta}.
\end{equation*}
Hence by (\ref{eq:u}) and (\ref{eq:v}) obtain
\begin{equation*}
 u_j(x)|x|^{2j}=\sum_{|\alpha|=m-2j}\sum_{|\beta|=j}\frac{j!}{\beta!}a_{j,\alpha}x^{\alpha+2\beta}
\end{equation*}
and
\begin{equation*}
v_k(x)|x|^{2k}=\sum_{|\tilde{\alpha}|=m-2k}\sum_{|\tilde{\beta}|=k}\frac{k!}{\tilde{\beta}!}b_{k,\tilde{\alpha}}x^{\tilde{\alpha}+2\tilde{\beta}}.
\end{equation*}
So, by Lemma \ref{Fischer} and (\ref{eq:Fischer}) we conclude that
\begin{multline*}
 \langle u, v \rangle_{F}=\sum_{j,k=0}^{p-1}\langle u_j(x)|x|^{2j}, v_k(x)|x|^{2k}\rangle_{F}\\
 =\sum_{j,k=0}^{p-1}\sum_{|\beta|=j}\sum_{|\tilde{\beta}|=k}
 \sum_{|\alpha|=m-2j \atop {|\tilde{\alpha}|=m-2k \atop \alpha+2\beta=\tilde{\alpha}+2\tilde{\beta}}}
 (\alpha+2\beta)!
 \frac{j!}{\beta!}\frac{k!}{\tilde{\beta}!}
 a_{j,\alpha}\overline{b_{k,\tilde{\alpha}}}.
\end{multline*}
\end{proof}

Directly from the above theorem we get
\begin{Wn}
 If $p>1$, then there does not exist a constant $C=C(n,m,p)$ depending only on $n,m$ and $p$ such that
 \begin{eqnarray*}
  \langle u, v \rangle_F= C \langle u, v \rangle_{\widehat{S}_p}\quad\textrm{for every}\quad
  u,v\in\mathcal{H}^p_m(\CC^n).
 \end{eqnarray*}
 In particular, if $u(x)=u_j(x)|x|^{2j}$ and $v(x)=v_j(x)|x|^{2j}$, then
 \begin{equation*}
  \langle u, v \rangle_{\widehat{S}_p}=\sum_{|\alpha|=m-2j}a_{j,\alpha}\overline{b_{j,\alpha}}\frac{\alpha!}{n(n+2)\cdots (n+2m-2(2j+1))}
 \end{equation*}
and 
\begin{equation*}
 \langle u, v \rangle_{F}=\sum_{|\alpha|=m-2j}\sum_{|\tilde{\alpha}|=m-2j}\sum_{|\beta|=j\atop {|\tilde{\beta}|=j \atop \alpha+2\beta=\tilde{\alpha}+2\tilde{\beta}}}
 (\alpha+2\beta)!
 \frac{j!}{\beta!}\frac{j!}{\tilde{\beta}!}
 a_{j,\alpha}\overline{b_{j,\tilde{\alpha}}}. 
\end{equation*}
\end{Wn}

\begin{Uw}
 In \cite{R} Render also considers the reproducing kernel $\tilde{Z}_m^p$ of $\mathcal{H}_m^p(\RR^n)$ with respect to the Fischer product (\ref{eq:Fischer}), which is defined by
 \begin{equation}
 \label{eq:Z.1}
\tilde{Z}_m^p(x,y):=\sum_{j=1}^{h_m^p}Q^j_m(x)\overline{Q^j_m(y)},
\end{equation}
where $\{Q^k_m(x)\}_{j=1}^{h_m^p}$ is an orthonormal basis of $\mathcal{H}_m^p(\RR^n)$ with respect to the Fischer product (\ref{eq:Fischer}). In particular he showed  the relation between the reproducing kernel $\tilde{Z}_m^p$ and zonal harmonics $Z_{m-2k}(x,y)$ (see  \cite[Theorem 2.4]{R}):
\begin{equation}
\label{eq:Z.2}
\tilde{Z}_m^p(x,y)=\sum_{k=0}^{\min\{[m/2],p-1\}}\frac{|x|^{2k}|y|^{2k}Z_{m-2k}(x,y)}{2^kk!n(n+2)\dots(n+2m-2k-2)}.
\end{equation}

Observe that $\tilde{Z}_m^p(x,y)$ one can treat as a version of zonal polyharmonic $Z_m^p(\eta,\zeta)$,
but with respect to the Fischer product (\ref{eq:Fischer}) instead of the integral product (\ref{inner product 1}). In particular (\ref{eq:Z.1})  and (\ref{eq:Z.2}) corresponds to properties of zonal polyharmonic $Z_m^p(\eta,\zeta)$ collected in Lemma \ref{St13}. Differences
between formulas given in Lemma \ref{St13} (ii) and in (\ref{eq:Z.2}) are consequences of the differences between the Fischer and integral products, which are described in Theorem \ref{Prop:difference}.
\end{Uw}

\section*{Acknowledgement}
The authors are grateful to the anonymous referees for the valuable comments and suggestions to improve the paper.

\end{document}